%% file: lifex-2.0-update.tex
\documentclass[preprint]{elsarticle}

\input{preamble}
\biboptions{sort&compress}

\begin{document}

\begin{frontmatter}
  \title{The lifex library version 2.0}

  \author[1]{Michele Bucelli\texorpdfstring{\corref{cor1}}{}}

  \affiliation[1]{
    organization={MOX Laboratory of Modeling and Scientific Computing, Dipartimento di Matematica, Politecnico di Milano},
    addressline={Piazza Leonardo da Vinci 32},
    postcode={20133},
    city={Milano},
    country={Italy}}

  \cortext[cor1]{Corresponding author. E-mail: michele.bucelli@polimi.it}
  \date{Last update: {\today}}

  \journal{}

  \begin{abstract}
    This article presents updates to \lifex{} [Africa, SoftwareX (2022)], a \cpp{} library for high-performance finite element simulations of multiphysics, multiscale and multidomain problems. In this release, we introduce an additional intergrid transfer method for non-matching multiphysics coupling on the same domain, significantly optimize nearest-neighbor point searches and interface coupling utilities, extend the support for 2D and mixed-dimensional problems, and provide improved facilities for input/output and simulation serialization and restart. These advancements also propagate to the previously released modules of \lifex{} specifically designed for cardiac modeling and simulation, namely \lifexfiber{} [Africa et al., BMC Bioinformatics (2023)], \lifexep{} [Africa et al., BMC Bioinformatics (2023)] and \lifexcfd{} [Africa et al., Computer Physics Communications (2024)]. The changes introduced in this release aim at consolidating \lifex{}'s position as a valuable and versatile tool for the simulation of multiphysics systems.
  \end{abstract}

\end{frontmatter}


\section{Introduction}

\begin{figure}
  \centering

  \includegraphics[width=0.5\textwidth]{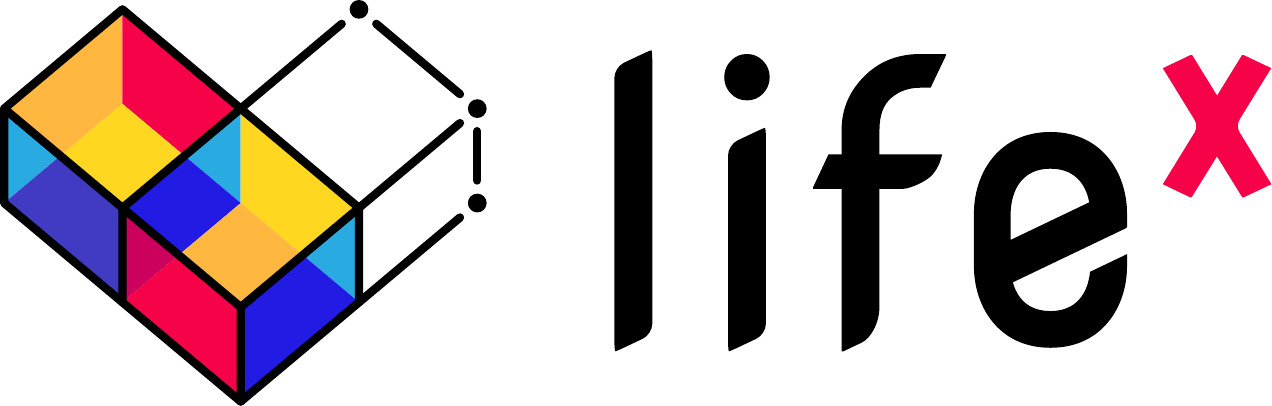}

  \caption{The \lifex{} logo. Image licensed under the CC-BY-SA 4.0 license.}
  \label{fig:logo}
\end{figure}

This paper discusses recent updates to \lifex{} \cite{africa2022lifexcore} (\url{https://lifex.gitlab.io}), a \cpp{} library tailored at finite element simulations for multiphysics and multiscale problems (logo depicted in \cref{fig:logo}). \lifex{} builds upon the finite element library \dealii{} \cite{arndt2019dealii,arndt2023dealii95} by implementing high-level reusable utilities for common tasks such as generating or reading meshes, both tetrahedral and hexahedral, solving linear and non-linear systems of equations, preconditioning linear systems, writing output to file, or managing simulation checkpointing and restart. Many of these features are implemented by wrapping \dealii{} functionality in easy-to-use software interfaces. These are configured through human-readable and well structured parameter files in the custom \dealii{} syntax \cite{africa2022lifexcore}, so that \lifex{} can be used effectively to design simulation tools that require minimal coding effort from the end user, if any \cite{africa2023lifexfiber,africa2023lifexep,africa2024lifexcfd}.

The effectiveness of the framework offered by \lifex{} was demonstrated by three modules specifically targeted at applications to cardiac and cardiovascular modeling: \lifexfiber{} \cite{africa2023lifexfiber}, for muscle fiber generation; \lifexep{} \cite{africa2023lifexep}, for simulating electrophysiology; \lifexcfd{} \cite{africa2024lifexcfd}, for \ac{CFD} simulations. The updates presented in this paper concern the general-purpose features of the library, and as such also apply to those modules.

\lifex{} provides the basis for a large and growing number of recent application studies \cite{barnafi2024reconstructing, bennati2023turbulent, bennati2024image, bucelli2023mathematical, capuano2024personalized, centofanti2024comparison, cicci2024efficient, corda2024influence, criseo2024computational, crispino2024cardiac, crugnola2024computational, crugnola2024inexact, de2024explicit, duca2024computational, falanga2024digital, fedele2023comprehensive, montino2024modeling, montino2024personalized, piersanti2024defining, renzi2024investigating, zappon2024integrated, zappon2024nonconforming, zingaro2023comprehensive, zingaro2024electromechanics, zingaro2024comprehensive}, which testify to the versatility and impact of the \lifex{} ecosystem on cardiovascular research. An up-to-date list of publications using \lifex{} is maintained on the dedicated page of the official website\footnote{\url{https://lifex.gitlab.io/lifex-public/publications.html}}.

The present release aims at consolidating, optimizing and enhancing the features of \lifex{}'s core module. We implement a new method for transferring data between non-matching meshes \cite{bucelli2023preserving, bucelli2024robust} in parallel simulations (\cref{sec:rbf}); we improve the interface and the performance of the parallel nearest-neighbor lookup of mesh vertices (\cref{sec:locators}), an operation used throughout the library for many tasks, including interface coupling; we introduce support for simulations in dimensions other than 3D (\cref{sec:dimensions}); we enhance simulation checkpointing and restart through a significantly improved and standardized user interface (\cref{sec:restart}), and introduce several improvements to the \ac{IO} facilities in general (\cref{sec:io}). All these features greatly enhance the usability and versatility of \lifex{}, with the aim of further improving its effectiveness as a tool for multiphysics simulations.

\Cref{tab:classes} reports the most relevant classes implemented in \lifex{}, grouped by category. The rest of this paper presents the features in the new release and their significance, while we refer to \cite{africa2022lifexcore} for a more extensive description of the general \lifex{} framework and of the features available in previous release, and to the source code repository\footnote{\url{https://gitlab.com/lifex/lifex-public}} and the online documentation\footnote{\url{https://lifex.gitlab.io/lifex-public/index.html}} for technical details. Unless otherwise specified, all \cpp{} classes discussed in this document are part of the \texttt{lifex::utils} namespace, which we omit henceforth in the interest of brevity.

\begin{table}
  \centering
  \footnotesize
  \ttfamily
  \begin{tabular}{l l l l}
    \toprule
    \textbf{Geometry}               & \textbf{Numerics}          & \textbf{Multiphysics}         & \textbf{IO} \\
    \midrule
    MeshHandler                     & NonLinearSolverHandler     & QuadratureEvaluation          & CSVWriter \\
    MeshInfo                        & LinearSolverHandler        & ProjectionL2                  & VTKFunction \\
    \textbf{DoFLocator} (*)         & PreconditionerHandler      & InterfaceHandler              & \textbf{OutputHandler} (*) \\
    \textbf{BoundaryDoFLocator} (*) & BlockPreconditionerHandler & \textbf{RBFInterpolation} (*) & \textbf{VTKImporter} (*) \\
    \textbf{GeodesicDistance} (*)   & BDFHandler                 &                               & \textbf{RestartHandler} (*) \\
                                    & BCHandler \\
                                    & FixedPointAcceleration \\
                                    & TimeInterpolation \\
                                    & \textbf{Laplace} (*) \\
    \bottomrule
  \end{tabular}

  \caption[]{List of \lifex{} core's most relevant classes, grouped by category. Classes that were introduced in the new release are written in bold and marked with an asterisk \texttt{(*)}. Additional details can be found in the online documentation (\url{https://lifex.gitlab.io/lifex-public/index.html}).}
  \label{tab:classes}
\end{table}

\section{Radial basis function interpolation for multiphysics coupling}
\label{sec:rbf}

\lifex{} is designed with multiphysics cardiac applications in mind, and as such has a strong focus on the coupling of heterogeneous models. Until previous release, the transfer of data between different domains or different models always imposed some conformity constraints on their discretization. Interface-coupled problems can be managed through the class \texttt{InterfaceHandler} and the related functionality, requiring that the two coupled models have a conforming discretization at their interface \cite{bucelli2023partitioned}. Data transfer between volume-coupled problems (that is, problems defined on the same domain) can be managed through the class \texttt{QuadratureEvaluation} and its derived classes, which require that the two models share the same mesh, although they can be discretized with different finite element spaces \cite{africa2022lifexcore}. The coupling of meshes of different resolution was, until now, only supported for nested grids of hexahedral elements \cite{regazzoni2022cardiac}.

With this release, we add support for \ac{RBF} interpolation between non-matching meshes, with the methods described in \cite{bucelli2023preserving, bucelli2024robust, deparis2014rescaled}. This allows to transfer data between spatial discretizations of arbitrary refinement, element shape (tetrahedral or hexahedral) and polynomial degree, thus greatly enhancing the library's flexibility in coupling heterogeneous problems. Furthermore, we support complex geometries with dedicated methods based on approximate geodesic distance \cite{bucelli2024robust}. This implementation of \ac{RBF} interpolation demonstrates excellent parallel scalability up to thousands of cores, as discussed in \cite{bucelli2023preserving, bucelli2024robust}.

The new features are exposed through the class \texttt{RBFInterpolation}, providing an interface for \ac{RBF} interpolation between arbitrary sets of points. Derived classes \texttt{RBFInterpolationDoFs} and \texttt{RBFInterpolationQuadrature} manage interpolation of data that is collocated at \acp{DoF} and mesh quadrature points, respectively. All these classes can be configured extensively through the parameter file. A new example named \texttt{ExampleRBFInterpolation} showcases the new features.

\section{Point locators}
\label{sec:locators}

\begin{figure}
  \centering

  \includegraphics{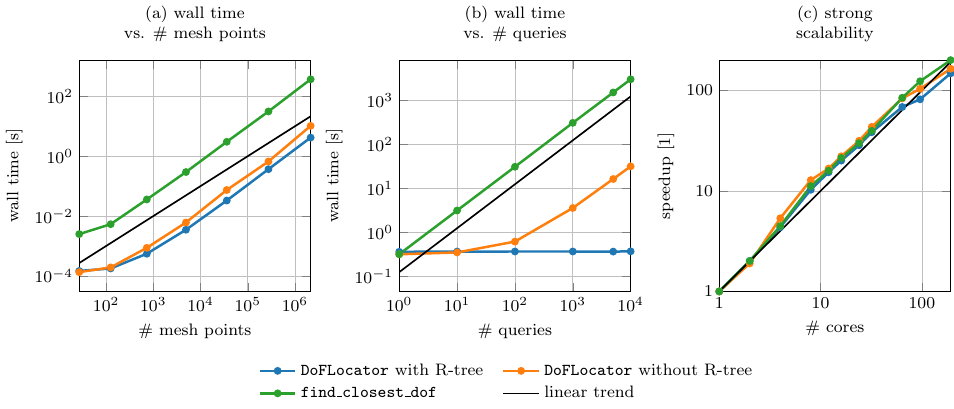}

  \caption{Performance of the new \texttt{DoFLocator} class compared to the old \texttt{find\_closest\_dof} function. (a) Wall time against the number of mesh points, for \num{100} repeated queries, using one parallel process. \texttt{DoFLocator}, both with and without R-trees, is close to \num{100} times faster than the old function. (b) Computational time against number of queries, for a mesh with \num{274625} vertices, using one parallel process. The advantage of R-trees becomes evident as the number of queries increases. (c) Nearly ideal parallel speedup of a strong scalability study, with \num{100} repeated queries on a mesh with \num{2146689} vertices.}
  \label{fig:locators}
\end{figure}

The release introduces new helper classes \texttt{DoFLocator} and \texttt{BoundaryDoFLocator} that offer an interface for the task of locating the nearest \ac{DoF} (or boundary \ac{DoF}) to a given point in the physical space, possibly in a parallel setting, which is a key part of many algorithms in \lifex{}.

The locator classes internally build an R-tree representation of the points \cite{manolopoulos2006rtrees}, using the implementation of R-trees from \texttt{boost::geometry::index} \cite{boost} wrapped by \dealii{}. Therefore, they provide a friendly yet computationally efficient interface for nearest-neighbor searches. A dedicated method is implemented for efficient multi-point queries, where every parallel process needs to locate a different set of points, possibly owned by other processes.

We stress that, until previous release, nearest-neighbor searches were done with a simple linear search algorithm, implemented by the \texttt{find\_closest\_dof} function (now removed). The current implementation significantly improves in terms of algorithmic complexity and performance, as shown in \cref{fig:locators}, and can take advantage from repeated queries by reusing the same R-tree.

\begin{figure}
  \centering

  \begin{subfigure}[b]{0.45\textwidth}
    \centering
    \includegraphics{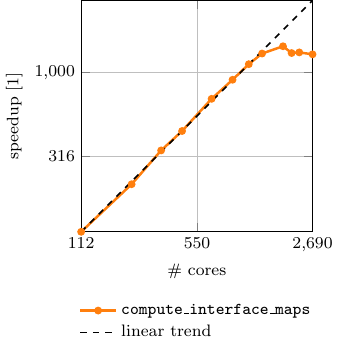}
    \caption{}
    \label{fig:interface-maps-scalability}
  \end{subfigure}
  \begin{subfigure}[b]{0.45\textwidth}
    \centering
    \includegraphics[width=0.95\textwidth]{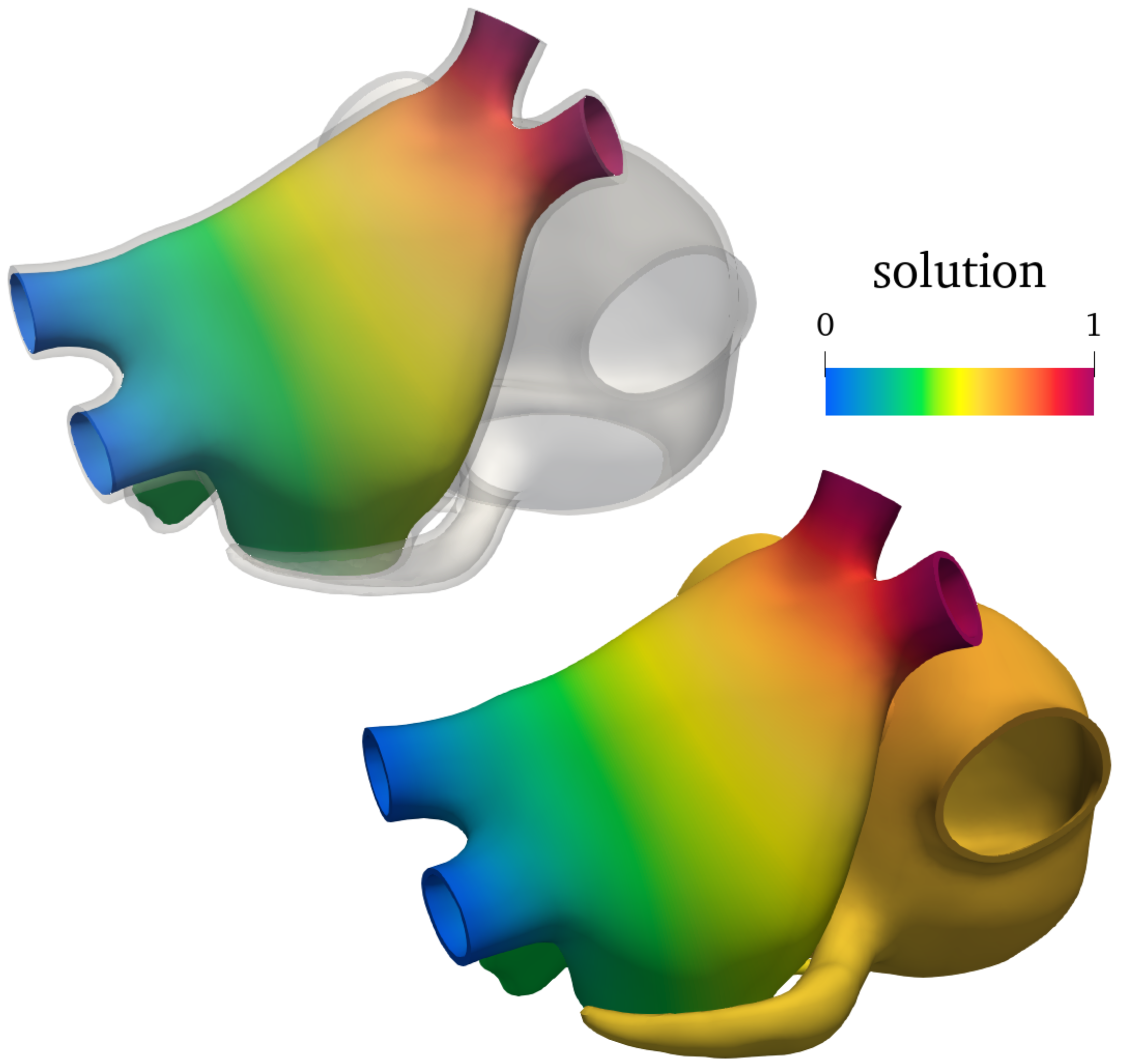}
    \caption{}
    \label{fig:2d-3d-coupling}
  \end{subfigure}

  \caption{(a) Strong scalability of the function \texttt{compute\_interface\_maps}, that establishes a mapping between the interface \acp{DoF} of two domains with a common boundary. This test couples a surface and a volume problem, with \num{1049601} and \num{135398529} \acp{DoF}, respectively, and \num{1049601} interface \acp{DoF}. (b) Example of coupling a surface and a volume problem through their common interface, as implemented by the new example \texttt{InterfaceCoupling2D3D}. We first solve a Laplace-Beltrami problem on a portion of the domain's boundary (top), and then solve a Laplace problem in the whole domain (bottom), with Dirichlet conditions taken from the surface problem. Atrial model taken from \cite{ferrer2015detailed}.}
\end{figure}

Additionally, the new locator classes are exploited in the coupling of domains across a common interface, as implemented by \texttt{InterfaceMap} and \texttt{InterfaceHandler}. The efficiency and parallel performance of the construction of interface maps has been significantly enhanced in this release. As depicted in \cref{fig:interface-maps-scalability}, the task of establishing a map between the interface \acp{DoF} of two domains with a common boundary (implemented by the function \texttt{compute\_interface\_maps}) shows ideal parallel scalability up to approximately \num{800} interface \acp{DoF} per process.

\section{1D, 2D and mixed-dimensional problems}
\label{sec:dimensions}

With this release, we improve \lifex{}'s support for 1D and 2D problems, and mixed-dimensional problems in general. The user can now specify the spatial dimension through the CMake \texttt{LIFEX\_DIM} parameter (which defaults to \num{3}). Many classes of general purpose are now templated over the physical and spatial dimensions, following the same convention as \dealii{}. These include \texttt{MeshHandler}, \texttt{MeshInfo}, \texttt{OutputHandlerBase} and its derivatives, \texttt{Laplace} and \texttt{GeodesicDistance}.

Most notably, with this improvement we extend the applicability of \lifex{} utilities to problems defined on surfaces (as seen in the new \texttt{LaplaceBeltramiExample}). Additionally, the \texttt{InterfaceHandler}-related utilities now allow to couple problems of mixed dimensions, such as a problem defined on a surface with another defined on a volume for which that surface is part of the boundary. This feature is demonstrated by the new example \texttt{InterfaceCoupling2D3D} (\cref{fig:2d-3d-coupling}).

\section{Checkpointing and restart}
\label{sec:restart}

\Ac{HPC} enviroments typically limit the maximum duration of a job to wall-times that are shorter than the duration of large-scale simulations. It is therefore crucial for a library such as \lifex{} to allow splitting computations over multiple jobs, overcoming these limitations.

To this end, \lifex{} allows to write the simulation state to a file (a process also referred to as \textit{checkpointing}, or \textit{serialization}), from which, at a later time, the simulation itself can be restarted. In this release, we significantly reworked this process, improving its robustness and reliability, enhancing its support for multiphysics simulations, and significantly simplifying its user interface.

All this is implemented through a new helper class named \texttt{RestartHandler}. The class can collect data from different problems, allowing to store multiple fields or scalar values in a single \texttt{.h5} file \cite{HDF5}. We remark that the previous release would create multiple files for each model, which could lead to confusion and clutter, and would not allow to include scalars in the serialized files, so that they would need to be saved and restored separately.

Conversely, the \texttt{RestartHandler} class offers a clean interface to store and retrieve all the data needed for restart. The class acts by keeping a list of references to data that needs to be serialized or deserialized. Such data can be easily registered through the methods \texttt{RestartHandler::attach\_scalar} and \texttt{RestartHandler::attach\_vector}, for scalar types and parallel vectors (or block vectors), respectively. A helper \texttt{RestartHandler::attach\_bdf\_handler} facilitates serialization and restart for time-dependent problems relying on the \texttt{BDFHandler} class \cite{africa2022lifexcore}. Serialization is done by calling the \texttt{RestartHandler::serialize} method, while restart is performed through the \texttt{RestartHandler::restart} method. Simulation and restart for multiple models can be easily centralized by having each model write to the same instance of \texttt{RestartHandler}.

All \lifex{} tutorials have been extended to exemplify the use of the \texttt{RestartHandler} class. Most notably, tutorial \num{6} demonstrates its use in the context of a multiphysics simulation in which the different sub-models are managed by separate classes.

From the user's perspective, checkpointing and restart are configured in two dedicated subsections of the parameter file:
\begin{lstlisting}[language=prm]
  subsection Serialization
    set Enable                      = true
    set Serialization basename      = restart
    set Serialize every n timesteps = 1000
  end

  subsection Restart
    set Enable                 = true
    set Restart basename       = out_dir/restart
    set Restart timestep index = 1000
  end
\end{lstlisting}
We stress that, differently from previous release, the user need not specify the initial time or initial timestep number of the restarting simulation, as these will be retrieved from the serialized data. Overall, this makes the process of restarting much simpler and less error-prone.

\section{Input/output enhancements}
\label{sec:io}

On top of the previously discussed \texttt{RestartHandler}, we introduced new classes to centralize \ac{IO} tasks that are common between multiple applications of \lifex{}. This has a twofold purpose: on the one hand, it enforces a standardized interface for those tasks, ensuring in particular that all applications share the same parameter file structure. On the other hand, this centralization greatly facilitates any future extension.

\subsection{Output of problem solutions to file}
Data output has been centralized to the new class \texttt{OutputHandler}, wrapping \dealii{}'s data writer class \texttt{dealii::DataOut}. With respect to previous release, we exposed output in \texttt{.pvtu/.vtu} format (on top of the already available \texttt{.xdmf/.h5} format). Indeed, we have observed that parallel output to \texttt{.h5} files may occasionally lead to deadlocks due to issues with parallel filesystems. The \texttt{.pvtu/vtu} format, where each process writes its data in an independent file, offers an effective workaround in those situations. We point out that \texttt{.pvtu/vtu} output usually occupies more disk space than \texttt{.xdmf/.h5} output, due to the latter allowing to filter out duplicate internal vertices.

\subsection{Reading data from VTK files}
Many applications are based on reading functional data to be used as parameters for numerical models implemented in \lifex{} \cite{capuano2024personalized, montino2024modeling, montino2024personalized}. A new class \texttt{VTKImporter} facilitates reading and remapping data from the well-established VTK file formats (\texttt{.vtk}, \texttt{.vtp} and \texttt{.vtu}), and optionally serializing the imported data to a binary file for later reuse. The class supports all the types of VTK functions offered by \texttt{VTKFunction} and \texttt{VTKPreprocess} \cite{africa2022lifexcore} (linear projections, closest-point projections and signed distance evaluation), but additionally takes care of standardizing the parameter file sections that configure these operations.

\subsection{Fixed memory occupation peak when reading meshes}
Until previous release, when reading tetrahedral meshes from file, all parallel processes would read the mesh in its entirety, and discard the portions attributed to other processes only after partitioning. This would lead to a very high peak memory occupation, often higher than the available memory, thus frequently resulting in the simulation being killed during initialization. We introduce a new parameter \texttt{Reading group size} to the \texttt{MeshHandler} class, which allows to reduce the number of processes that read the entire mesh (based on \texttt{create\_description\_from\_triangulation\_in\_groups} from \dealii{}'s \texttt{TriangulationDescription::Utilities} namespace). This has proven crucial in supporting very large-scale simulations.

\section{Additional improvements}

In addition to multiple bufgixes and performance improvements, the new release includes the following changes:
\begin{itemize}
  \item \lifex{} is now updated to use \dealii{} version 9.5.1;

  \item a new helper class \texttt{Laplace} provides a simple interface for solving Laplace and Laplace-Beltrami problems (that is, $-\Delta u = 0$ in a certain domain $\Omega$). The class is meant to be used for algorithms that require solving the Laplace equation as an intermediate step, such as the ones discussed in \cite{africa2023lifexfiber};
  
  \item users can specify a custom set of default parameters to each instance of \texttt{PreconditionerHandler}. This is particularly useful since the optimal preconditioner configuration may vary between different problems: such a configuration can be built into the source code for each problem, without requiring the user to manually adjust the parameter file;

  \item a new class \texttt{GeodesicDistance} allows to compute an edge-based approximation of the geodesic distance and of the shortest path between points within a domain, exposed through the new app \texttt{shortest\_path}.
\end{itemize}

\section{Conclusions}
\label{sec:conclusions}

The \lifex{} library offers a comprehensive set of tools to facilitate the development of multiphysics finite element simulations. The 2.0 release described in this paper extends the library's functionality, by improving its multiphysics coupling capabilities and its support for simulations of different dimensionalities. Furthermore, it improves the efficiency and parallel scalability of fundamental tasks such as parallel point location and communication between interface-coupled problems. Finally, the new release provides an improved and standardized interface for several tasks related to input and output. All these changes significantly enhance the capabilities of \lifex{} and its applicability to large-scale problems, and consolidate its position as a valuable framework for simulating multiphysics systems.

\section*{Acknowledgements}

\lifex{} was developed under the scientific supervision of Profs. Luca Dede' (Politecnico di Milano, Milano, Italy) and Alfio Quarteroni (Politecnico di Milano, Milano, Italy). The author wishes to thank all the \lifex{} users and developers for their software contributions, testing and feedback, and especially Pasquale Claudio Africa (SISSA, Trieste, Italy) for his years as \lifex{} maintainer and his invaluable support. The \lifex{} logo was designed by Silvia Pozzi (Artiness).

The present research is part of the activities of ``Dipartimento di Eccellenza 2023--2027'', MUR, Italy, Dipartimento di Matematica, Politecnico di Milano. The author has received support from the project PRIN2022, MUR, Italy, 2023--2025, 202232A8AN ``Computational modeling of the heart: from efficient numerical solvers to cardiac digital twins''. The author acknowledges his membership to INdAM GNCS - Gruppo Nazionale per il Calcolo Scientifico (National Group for Scientific Computing, Italy), and INdAM GNCS project CUP E53C23001670001. The author acknowledges ISCRA for awarding this project access to the LEONARDO supercomputer, owned by the EuroHPC Joint Undertaking, hosted by CINECA (Italy).

\bibliographystyle{abbrv}
\bibliography{bibliography}

\end{document}

%% file: preamble.tex
\usepackage[english]{babel}
\usepackage[htt]{hyphenat}
\usepackage{acro}

\usepackage[a4paper, total={6.5in, 8.66in}]{geometry}
\usepackage{subcaption}

\usepackage{booktabs}

\usepackage{amsmath}
\usepackage{amssymb}
\usepackage{mathtools}
\usepackage{amsthm}
\usepackage{siunitx}

\usepackage[hypertexnames=false,
            colorlinks=true,
            pdfborderstyle={/S/U/W 0},
            citecolor=blue]{hyperref}
\usepackage[capitalise, noabbrev]{cleveref}

\usepackage{listings}

\usepackage[colorinlistoftodos,
            bordercolor=white,
            textsize=footnotesize]{todonotes}
\usepackage{xcolor}

\definecolor{code-bg}{rgb}{0.93, 0.93, 0.97}
\definecolor{code-comment}{rgb}{0.5, 0.5, 0.5}
\definecolor{code-keyword}{rgb}{0.75, 0.0, 0.25}
\definecolor{plot1}{HTML}{1f77b4}
\definecolor{plot2}{HTML}{ff7f0e}
\definecolor{plot3}{HTML}{2ca02c}
\definecolor{plot4}{HTML}{d62728}
\definecolor{plot5}{HTML}{9467bd}

\lstdefinestyle{codestyle}{
    basicstyle=\ttfamily\footnotesize,
    backgroundcolor=\color{white},
    commentstyle=\color{code-comment},
    keywordstyle=\color{code-keyword},
    breakatwhitespace=false,
    breaklines=true,
    captionpos=b,
    keepspaces=true,
    numbersep=5pt,
    showspaces=false,
    showstringspaces=false,
    showtabs=false,
    tabsize=2
}
\lstdefinelanguage{prm}{
    keywords={subsection, set, end},
    comment=[l]{\#}
}

\newcommand{\lifex}{\texttt{life\textsuperscript{x}}}
\newcommand{\lifexfiber}{\texttt{life\textsuperscript{x}-fiber}}
\newcommand{\lifexep}{\texttt{life\textsuperscript{x}-ep}}
\newcommand{\lifexcfd}{\texttt{life\textsuperscript{x}-cfd}}
\newcommand{\dealii}{\texttt{deal.II}}
\newcommand{\cpp}{\texttt{C++}}

\DeclareAcronym{CFD}{short={CFD}, long={computational fluid dynamics}}
\DeclareAcronym{IO}{short={IO}, long={input/output}}
\DeclareAcronym{RBF}{short={RBF}, long={radial basis function}}
\DeclareAcronym{DoF}{short={DoF}, long={degree of freedom}, long-plural-form={degrees of freedom}}
\DeclareAcronym{HPC}{short={HPC}, long={high-performance computing}}

\makeatletter
\def\ps@pprintTitle{%
 \let\@oddhead\@empty
 \let\@evenhead\@empty
 \def\@oddfoot{}%
 \let\@evenfoot\@oddfoot}
\makeatother